\documentclass[12pt]%
{amsart}  
\usepackage{verbatim}
\usepackage{amsmath} 
\usepackage{amsthm}
\usepackage{amssymb}
\usepackage{enumerate}
\newtheorem{theorem}{Theorem}[section]

\newtheorem{lemma}[theorem]{Lemma}

\newtheorem{claim}{Claim}[theorem]

\newtheorem{obse}[theorem]{Observation}
\newtheorem{qclaim}{Claim}  

\theoremstyle{definition}

\newtheorem{definition}[theorem]{Definition}

\theoremstyle{remark}
  
{}

\newcommand{\eproof}{\end{pf}}
\newcommand{\Eproof}{\end{pf*}}
\newcommand{\sproof}[1]{\begin{pf*}{#1}}
\newcommand{\esproof}{\end{pf*}}

\def\myheads#1;#2;{
\pagestyle{myheadings}
\markboth{{\sc\hfill #1\hfill\protect\makebox[0cm][r]{\rm\today}}}
{{\sc\protect\makebox[0cm][l]{\rm\today}\hfill #2\hfill}}
}

\newcommand{\acal}{{\mathcal A}}
\newcommand{\bcal}{{\mathcal B}}
\newcommand{\ccal}{{\mathcal C}}
\newcommand{\dcal}{{\mathcal D}}

\newcommand{\fcal}{{\mathcal F}}

\newcommand{\hcal}{{\mathcal H}}

\newcommand{\mcal}{{\mathcal M}}

\newcommand{\setm}{\setminus}
\newcommand{\empt}{\emptyset}
\newcommand{\subs}{\subset}

\newcommand{\oo}{{{\omega}_1}}
\newcommand{\dom}{\operatorname{dom}}
\newcommand{\ran}{\operatorname{ran}}
\def\<{\left\langle}
\def\>{\right\rangle}

\def\OO{{\omega}}
\def\oo{\omega_1}
\def\br#1;#2;{\bigl[ {#1} \bigr]^ {#2} }
\def\bc#1;#2;{\bigl( {#1} \bigr)^ {#2} }

\def\ooseq#1;#2;{\< {#1}_{#2}:{#2}<\oo\>}
\def\ooset#1;#2;{\{ {#1}_{#2}:{#2}<\oo\}}
\def\seq#1;#2;#3;{\< {#1}_{#2}:{#2}<#3\>}
\def\set#1;#2;#3;{\{ {#1}_{#2}:{#2}<#3\}}
\def\oseq#1;#2;{\< {#1}_{#2}:{#2}<\OO\>}
\def\oset#1;#2;{\{ {#1}_{#2}:{#2}<\OO\}}
\def\oosequ#1;#2;{\< {#1}^{#2}:{#2}<\oo\>}
\def\oosetu#1;#2;{\{ {#1}^{#2}:{#2}<\oo\}}
\def\sequ#1;#2;#3;{\< {#1}^{#2}:{#2}<#3\>}
\def\setu#1;#2;#3;{\{ {#1}^{#2}:{#2}<#3\}}
\def\osequ#1;#2;{\< {#1}^{#2}:{#2}<\OO\>}
\def\osetu#1;#2;{\{ {#1}^{#2}:{#2}<\OO\}}

\def\force{\raisebox{1.5pt}{\mbox {$\scriptscriptstyle\|$}}
\mbox{$\!\mbox{---}$}}

\newcommand{\fn}{\operatorname{Fn}}

\def\to{\longrightarrow}

\newcommand{\id}{\operatorname{id}}

\def\fin#1;{\br #1;<{\omega};}

\newcommand{\restr}%
{\mathop{\hspace{0.01ex}|\hspace*{-0.02ex}{\grave{}}\hspace{0.4ex}}}

\newcommand{\rest}{\restr}

\newcommand{\prtime}{{\count0=\time\divide\count0 by 60
\count1=-\count0\multiply\count1 by 60
\advance\count1 by \time
\the\count0:\the\count1}
}

\def\myheads#1;#2;{
\pagestyle{myheadings}
\markboth{{\sc\hfill #1\hfill\protect\makebox[0cm][r]{\rm\today; \prtime}}}
{{\sc\protect\makebox[0cm][l]{\rm\today;\ \prtime}\hfill #2\hfill}}
\thispagestyle{myheadings}
}

\let\QED\qed
\newcommand{\prlabel}[1]{\renewcommand{\qed}{\QED${}_{\ref{#1}}$}}
\newcommand{\prtxtlabel}[1]{\renewcommand{\qed}{\QED${}_{{\mbox{\tiny #1}}}$}}
\newcommand{\prnolabel}{\prtxtlabel{}}

\newcommand{\ffcal}[1]{\fcal(#1)}
\newcommand{\sscal}[1]{\fn(#1,{\Omega})}
\newcommand{\mmem}[3]{#1_{{#2,#3}}}
\newcommand{\mdmem}[3]{\mmem{\dot{#1}}{#2}{#3}}

\newcommand{\amem}[2]{\mmem {A}{#1}{#2}}
\newcommand{\dmem}[2]{\mmem {D}{#1}{#2}}

\newcommand{\nfam}[1]{N(#1)}
\newcommand{\kfun}[1]{K^{#1}}
\newcommand{\ddmem}[2]{\mdmem {D}{#1}{#2}}

\newcommand{\oot}{{\omega}_2}
\newcommand{\nice}{nicely extendible}                                       
\newcommand{\good}{${\omega}$-determined}                                       

\newcommand{\blo}[1]{\mathbb I_{#1}}

\newcommand{\matr}[2]{\mcal(#1,#2)}
\newcommand{\omatr}[1]{\matr {#1}{\omega}}
\newcommand{\pto}{\stackrel{p}{\to}}
\newcommand{\ma}{\operatorname{{\gamma}}}
\newcommand{\coh}[1]{\ccal_{#1}}
\newcommand{\siso}[2]{{\sigma}_{#1,#2}}

\newcommand{\pin}[1]{\operatorname{\pi_0}(#1)}
\newcommand{\pie}[1]{\operatorname{\pi_1}(#1)}
\newcommand{\pii}[1]{\operatorname{\pi_i}(#1)}
\newcommand{\pip}[2]{\operatorname{\pi_{#2}}(#1)}

\def\lev#1;#2;{\operatorname{I}_{#1}(#2)}
\newcommand{\htt}{\operatorname{ht}}
\newcommand{\wdt}{\operatorname{wd}}

\begin{document}

\title{Cardinal sequences and Cohen real extensions}

\author[I. Juh\'asz]{Istv\'an Juh\'asz}
\thanks{The first, third and fourth authors were supported by the 
Hungarian National Foundation for Scientific Research grant no. 37758
}   
\address{Alfr{\'e}d R{\'e}nyi Institute of Mathematics}
\email{juhasz@renyi.hu}

\author[S. Shelah]{Saharon Shelah}
\thanks{The second author was supported by The Israel Science
  Foundation founded by the Israel Academy of Sciences and Humanities. 
  Publication 765}
\address{The Hebrew University of Jerusalem}
\email{shelah@math.huji.ac.il}

\author[L. Soukup]{
Lajos Soukup }
\thanks{
The third author was partially
supported by Grant-in-Aid for JSPS Fellows No.\ 98259 of the Ministry of 
Education, Science, Sports and Culture, Japan}
\address{Alfr{\'e}d R{\'e}nyi Institute of Mathematics }  
\email{soukup@renyi.hu}

\author[Z. Szentmikl\'ossy]{
Zolt\'an Szentmikl\'ossy}
\address{E\"otv\"os University of Budapest} 
\email{zoli@renyi.hu}

\begin{abstract}
We show that if we add any number of Cohen reals to the ground model then, in
the generic extension,  a locally compact
scattered space has at most 
$(2^{\aleph_0})^V$
many levels of size $\omega$.

We also give a complete $ZFC$  characterization of the cardinal sequences of 
regular scattered spaces. Although the classes of the regular and 
of the $0$-dimensional scattered spaces are different, 
we prove that they
have the same cardinal sequences.
\end{abstract}

\maketitle
\section{Introduction}\label{sc:SH}

Let us start by recalling that a topological space $X$ is called 
{\em scattered} if every non-empty 
subspace of $X$ has an isolated point. Via the well-known Cantor-Bendixson
analysis then $X$ decomposes into levels,  
the ${\alpha}^{\text{th}}$ 
Cantor-Bendixson level of  $X$ will be 
denoted by $\lev {\alpha};X;$. 
The {\em height of $X$, $\htt(X)$,}
is the least ordinal ${\alpha}$ with $\lev {\alpha};X;=\empt$.
The {\em width of $X$, $\wdt(X)$,}
is defined by  
$\wdt(X) = \sup\{\ |\lev {\alpha};X;|:{\alpha}<\htt(X)\}$.
Our main object of study is
the {\em cardinal sequence} of  $X$, 
denoted by $\operatorname{CS}(X)$, that is the
sequence of cardinalities of the
non-empty Candor-Bendixson levels of $X$, i.e.
\begin{displaymath}
CS(X)=\bigl\langle\ |I_{\alpha}(X)|:{\alpha}<\htt(X)\ \bigr\rangle.
\end{displaymath}

The cardinality of a $T_3$ , in particular of a 
locally compact,  scattered  $T_2$ (in short: LCS) space $X$
is at most $2^{\, |\lev 0;X;| }$, hence clearly 
$\htt(X)<(2^{|\lev 0;X;|})^+$
and $|I_{\alpha}(X)| \le 2^{|\lev 0;X;|}$ for each $\alpha$. 
(Locally compact scattered spaces are
closely related to superatomic boolean algebras via Stone duality and
the study of their cardinal sequences was actually originated in
that subject.) 
Thus, in particular,  under $CH$ there is no
scattered $T_3$ space of height ${\omega}_2$ 
and having only countably many isolated points.    
After I. Juh\'asz and W.  Weiss, 
\cite[theorem 4]{JW}, had proved in ZFC that
for every ${\alpha}<{\omega}_2$ there is an 
LCS space $X$ with $\htt(X)={\alpha}$ and $\wdt(X)={\omega}$,
it was a  natural question if the existence of
an LCS space of height ${\omega}_2$ and  width ${\omega}$
follows from $\neg CH$.
 This question
was answered in the negative  by W. Just who proved,
\cite[theorem 2.13]{J1}, that if one blows up the continuum 
by adding Cohen reals 
to a model of $CH$ then in the resulting generic extension there is no
LCS space of height ${\omega}_2$
and width ${\omega}$.  
On the other hand, in their ground breaking work \cite{BS},   
J. Baumgartner and
S. Shelah produced a model in which there is a LCS space of height
${\omega}_2$ and width ${\omega}$, moreover they 
proved  in ZFC that for each ${\alpha}<(2^{\omega})^+$
there is a scattered 0-dimensional $T_2$ space $X$
with $\htt(X)={\alpha}$ and $\wdt(X)={\omega}$. 
Building on the idea of the proof of this latter
result, in  section \ref{sc:regular} 
we succeeded in giving  
a complete characterization of the cardinal sequences of both 
$T_3$ and zero-dimensional $T_2$ scattered spaces. 
Although the classes of the regular and of the zero-dimensional
scattered spaces are different, it will turn out that they
yield the same class of cardinal sequences.
We should add that, with quite a bit of extra effort,
in \cite{M}, J.-C. Martinez
extended the former result of Baumgartner and Shelah by producing a
model in which for every ordinal $\alpha < \omega_3$ there is a LCS
space of height $\alpha$ and width $\omega$. The question if it is
consistent to have a LCS space of height $\omega_3$ and width $\omega$
remains a big mystery.

In section \ref{sc:cohen} we strengthened the result of Just by proving,
in particular, that
in the same Cohen real extension no LCS space may 
have $\oot$ many countable (non-empty) levels. It seems to be an intriguing
(and natural) problem if the non-existence of an LCS space of width $\omega$
and height $\oot$ implies in ZFC the above conclusion, or more generally:
is any subsequence of the cardinal sequence of an LCS space again such
a cardinal sequence?
In connection with this problem let us remark that,
(as is shown in \cite{Fu} or \cite{JK}), in 
the side-by-side random real extension of a model of CH  
the combinatorial principle $\ccal^s({\oot})$
introduced in \cite[definition 2.3]{JSSz}
holds, consequently in such an extension
there is no LCS space $X$ of height ${\omega}_2$
and width ${\omega}$. In fact,  by \cite[theorem 4.12]{JSSz}, 
 $\ccal^s({\oot})$ implies that 
$\{\alpha \in \omega_2 : |I_\alpha(X)| = \omega\}$ is non-stationary
in ${\omega}_2$.
However,
we do not know  if our above mentioned result, namely 
theorem \ref{tm:scat}, holds there.

The morale of our above discussion may be concisely formulated as follows:
The cardinal sequences of regular or zero-dimensional 
scattered spaces are only
subject to the trivial inequality $|X| \le 2^{|I_0(X)|}$, 
however those of the LCS spaces are much harder to determine, in particular,
they are sensitive to the model of set theory in 
which we look at them.


\section{Countable levels in Cohen real extensions}\label{sc:cohen}

Let us formulate then the promised strengthening of Just's result. We
note that no assumption (like CH) is made on our ground model.  

\begin{theorem}\label{tm:scat}
Let us set ${\kappa}=(2^{\omega})^+$  
and add any number of Cohen reals to our ground model. Then in
the resulting extension no LCS space contains a $\kappa$-sequence
$\{E_{\alpha}:{\alpha}<{\kappa}\}$ of  
pairwise disjoint countable  subspaces such that  
$\overline{E_{\alpha}}\supset E_{\beta}$ holds for all 
${\alpha}<{\beta}<{\kappa}$.
In particular, for any LCS space $X$ we have 
$\bigl|\{{\alpha}:|\lev {\alpha};X;|={\omega}\}\bigr|<{\kappa}$.
\end{theorem}

In fact, we shall prove a more general statement, but to formulate that
we need a definition.
A family of pairs (of sets)
$\dcal=\bigl\{\<D_0^{\alpha},D_1^{\alpha}\>:{\alpha}\in I
\bigl\}$ 
is said   to be {\em dyadic over a set $T$} iff
$D^{\alpha}_0\cap D^{\alpha}_1=\empt$ for each ${\alpha}\in I$
and 
$$
\dcal[\varepsilon]=
\bigcap\{D_{{\varepsilon}({\alpha})}^{\alpha}:{\alpha}\in\dom
{\varepsilon} \}$$ intersects $T$
for each ${\varepsilon}\in \fn(I,2)$.
We simply say that {\em $\dcal$ is dyadic} iff it is dyadic for some
$T$, i.e. $\dcal[{\varepsilon}]\ne \empt$ for each 
${\varepsilon}\in \fn(I,2)$
.

Now, it is obvious that in a LCS space \begin{itemize}
\item the compact open sets form a
base that is closed under finite unions,
\item there is no
infinite dyadic system of pairs of compact  sets.
\end{itemize}
Consequently, theorem \ref{tm:dia} below
immediately yields theorem \ref{tm:scat} above.


\begin{theorem}\label{tm:dia}
Set ${\kappa}=(2^{\omega})^+$  
and add any number of Cohen reals to the ground model. Then in
the resulting generic extension  the following statement holds: 
If $X$ is any $T_2$ space 
 containing  pairwise disjoint countable subspaces  
$\{E_{\alpha}:{\alpha}<{\kappa}\}$ such that 
$\overline{E_{\alpha}}\supset E_{\beta}$ for 
${\alpha}<{\beta}<{\kappa}$ and $X=\overline {E_0}$ (i. e. ${E_0}$ is dense
in $X$) , 
moreover,  for each $x\in X$, we have fixed a  
neighbourhood base $\bcal(x)$
of $x$ in $X$  that is closed under finite unions then there is an infinite
set $a\in \br {\kappa};{\omega};$,  for
each ${\alpha}\in a$ there are disjoint finite subsets 
$L^0_{\alpha}$ and $L^1_{\alpha}$ of
$E_{\alpha}$, and for each $x\in L^0_{\alpha}\cup L^1_{\alpha}$ there is a
basic neighbourhood $V(x) \in \bcal(x)$
such that the infinite family of pairs
\begin{equation}\notag
\Bigl\{\bigl\langle
\bigcup_{x\in L^0_{\alpha}}V(x),\bigcup_{x\in L^1_{\alpha}}V(x)
\bigr\rangle:{\alpha}\in a\Bigr\}
\end{equation}
is dyadic.
\end{theorem}

This topological statement in the Cohen extension
in turn will follow   from a purely combinatorial one concerning
certain matrices, namely theorem \ref{tm:mc}.

To formulate this theorem we again
need  some notation and definitions.

For an ordinal ${\alpha}$ 
the interval $[{\omega}{\alpha},{\omega}{\alpha}+{\omega})$
will be denoted by  $\blo {\alpha}$.

Given two sets  $A$ and $B$
we write $f:A\pto B$ to denote that 
$f$ is a partial  function from $A$ to $B$, i. e.
a function from a subset of $A$ into $B$.
As usual, we let 
\begin{displaymath}
\fn(A,B)=\{f : |f|<{\omega} \text{ and } f:A\pto B\}.
\end{displaymath}

If $A\subs{\rm On}$ then
for any partial function $f:A\pto B$ we set
\begin{displaymath}
\ma (f)=\left \{ 
\begin{array}{ll}
\min\dom f&\text{if $\dom f\ne\empt$,}\\
\sup A&\text{if $\dom f=\empt$.}
\end{array}
\right .
\end{displaymath}

We let
\begin{equation}\notag
\Omega=\bigl\{\<A,B\>\in \br{\omega};<{\omega};\times\br {\omega};<{\omega};:
A\cap B=\empt \bigr\},
\end{equation}
and for $\ell=\<A,B\>\in \Omega$ we set $\pin \ell=A$ and $\pie \ell=B$.

If $S$ and $T$ are  sets of ordinals, we denote by 
$\matr {S}{T}$ the family of all $S\times{\omega}$-matrices consisting 
of subsets of $T$, i. e. $\acal\in\matr ST$
means that $\acal=\<\amem {\alpha}i:{\alpha}\in S
,i\in{\omega}\>$,
where $\amem {\alpha}i\subs T$ for each ${\alpha}\in S$ 
and $i<{\omega}$. 



For $\acal\in \matr {S}{T}$, $f:S\pto S$, and  $s:S\pto {\Omega}$ 
the pair $(f,s)$ is said to be {\em $\acal$-dyadic  (over  $U$)} iff the family
of pairs
\begin{multline}\notag
\Bigl\{\Bigl\langle
\cup\bigl\{{\amem{f({\alpha})}{n}:n\in \pin {s({\alpha})}}\bigr\},
\cup\bigl\{{\amem{f({\alpha})}{n}:n\in \pie {s({\alpha})}}\bigr\}
\Bigr\rangle:\\{\alpha}\in\dom f\cap \dom s\Bigr\}.
\end{multline}
is dyadic (over $U$).
If the pair 
$\<\id_S,s\>$ is $\acal$-dyadic (over $U$)
then    $s$ is simply called {\em $\acal$-dyadic (over  $U$)}.
It is this latter notion of $\acal$-dyadicity of a single partial function
that is really important (that for pairs is only of technical significance).
Hence we state below an alternative characterisation of it.

For $\acal\in\matr ST$, 
$s:S\pto\Omega$, and  
${\varepsilon}\in \fn(\dom s,2)$ we write
\begin{displaymath}\notag
\acal[s,{\varepsilon}]=\bigcap_{{\alpha}\in \dom {\varepsilon}}
\bigcup\{\amem {{\alpha}}n:n\in \pip {s({\alpha})}{{\varepsilon}({\alpha})}\}.
\end{displaymath}

\begin{obse}
If $\acal\in\matr ST$ then
$s:S\pto\Omega$ is $\acal$-dyadic over $U$ 
iff $\acal[s,{\varepsilon}]\cap U\ne\empt$
for each ${\varepsilon}\in \fn(\dom s,2)$ and 
\begin{displaymath}
\bigcup\{\amem{\alpha}n:n\in \pin{s({\alpha})}\}\cap
\bigcup\{\amem{\alpha}n:n\in \pie{s({\alpha})}\}
=\empt
\end{displaymath}
for each ${\alpha}\in \dom s$.
\end{obse}

The following easy observation will be applied later, in the proof of
lemma \ref{lm:nice}:
\begin{obse}\label{obse:circ}
If $g:S\pto S$ and $s:S\pto {\Omega}$ satisfy $\dom s\subs \ran g$, and the
pair $(g,s\circ g)$ is $\acal$-dyadic over $U$ then  $s$
 is $\acal$-dyadic over $U$, as well.
\end{obse}

\begin{definition}\label{df:md}
Fix a cardinal ${\kappa}$ and
let $\dcal\in \matr{\kappa}{\kappa}$.
For $s:{\kappa}\pto {\Omega}$
we say that ${s}$ is {\em $\dcal$-min-dyadic (m.d.)} iff
$s$ is $\dcal$-dyadic over $\blo {\ma(s)}$.

Moreover, we say that the matrix
$\dcal$ is {\em m.d.-extendible } iff for each finite $\dcal$-min-dyadic 
partial function $s:{\kappa}\pto {\Omega}$
and  for each ${\gamma}<\ma(s)$
there is an $\ell\in {\Omega}$ such that $s\cup\{\<{\gamma},\ell\>\}$
is also $\dcal$-min-dyadic, i. e. $\dcal$-dyadic over $\blo {\gamma}$ .
\end{definition}

Since $\blo 0={\omega}$, we clearly have the following.
\begin{obse}\label{lm:mdd}
If $\dcal\in \matr{\kappa}{\kappa}$ is m.d-extendible 
and $s:{\kappa}\pto {\Omega}$ is a finite $\dcal$-min-dyadic 
partial function then
$s$ is $\dcal$-dyadic over  ${\omega}$.
\end{obse}


Finally, a  matrix $\dcal\in \matr {\kappa}{\kappa}$ 
will be called {\em\good} iff
$\dmem {\alpha}n\cap \dmem {\alpha}m\cap {\omega}=\empt$
implies $\dmem {\alpha}n\cap \dmem {\alpha}m=\empt$
whenever ${\alpha}<{\kappa}$ and $n<m<{\omega}$.

With this we now have all the necessary ingredients to formulate and 
prove the promised combinatorial statement that will be valid in any
Cohen real extension.

\begin{theorem}\label{tm:mc}
Set ${\kappa}=(2^{\omega})^+$    
and add any number of Cohen reals to the ground model. Then in
the resulting generic extension for every \good\ and m.d.-extendible matrix  
$\dcal\in\matr {\kappa}{\kappa}$ there is  an infinite $\dcal$-dyadic
partial  function  $h:{\kappa}\pto {\Omega}$.
\end{theorem}

Before proving theorem \ref{tm:mc}, however, we show how theorem
\ref{tm:dia} can be deduced from it.

\begin{proof}[Proof of theorem \ref{tm:dia} using theorem \ref{tm:mc}]
\prlabel{tm:dia}

We can assume without any loss of generality that  
$E_{\alpha}=\blo {\alpha}$  for each
${\alpha}<{\kappa}$ and then will define an appropriate matrix
$\dcal\in \matr {\kappa}{\kappa}$.

To this end, for coding purposes,  we first fix a bijection ${\rho}:
\br {\omega};2;\to {\omega}$ and let
 ${\eta}:{\omega}\to {\omega}$ and 
${\nu}:{\omega}\to {\omega}$ be the "co-ordinate" functions 
of its inverse, i. e.  
 $k={\rho}(\{{\nu}(k),{\eta}(k)\})$ and ${\nu}(k)<{\eta}(k)$ for each 
$k<{\omega}$.

Since $X$ is $T_2$, for each $n<{\omega}$  we can simultaneously pick 
basic  neighbourhoods $B^{\alpha}_n(m)\in\bcal({\omega}{\alpha}+m)$  of the
  points ${\omega}\cdot{\alpha}+m\in E_{\alpha}=\blo {\alpha}$ for
all $m< n$ such that the  sets  $\{B^{\alpha}_n(m)$: $m<n\}$ 
are pairwise disjoint.
 
Now we define
$\dcal=\<\dmem{{\alpha}}{k}:\<{\alpha},k\>\in {\kappa}\times
{\omega}\>\in \matr {\kappa}{\kappa}$ as follows:
\begin{equation}\notag
\dmem{\alpha}k=
B^{\alpha}_{{\eta}(k)}({\nu}(k))\cap {\kappa}. 
\end{equation}

This matrix $\dcal$ is clearly \good\ because
$E_0=\blo 0={\omega}$ is dense in $X$. It is a bit less easy to establish
the following

\begin{qclaim}\label{lm:dext}
$\dcal$ is also m.d.-extendible.
\end{qclaim}

\begin{proof}[Proof of the claim]\prlabel{lm:dext}
Let   $s:{\kappa}\pto {\Omega}$ be a finite $\dcal$-min-dyadic partial
function and 
let ${\gamma}<\ma(s)$.

Since the sets
$\{\dcal[s,{{\varepsilon}}]:{\varepsilon}\in {}^{\dom s}2\}$ are all open in
the subspace 
${\kappa}$ and they all intersect $\blo {\ma(s)}$, moreover 
every element of $\blo{\ma(s)}$ is an accumulation point of ${\blo
{\gamma}}$, it  follows that 
$\dcal[s,{{\varepsilon}}]\cap \blo {\gamma}$ must be infinite for each 
${\varepsilon}\in {}^{\dom s}2$.
Thus we can easily pick two disjoint finite subsets $A_0$ and $A_1$ of $\blo{\gamma}$
such that every $\dcal[s,{\varepsilon}]$ intersects both $A_0$ and $A_1$.
Let $n<{\omega}$ be chosen in such a way that 
$A_0\cup A_1\subs \{{\omega}{\gamma}+m:m<n\}$, and set
 $K_i=\{{\rho}\{m,n\}:m<n\land {\omega}{\gamma}+m\in A_i\}$ for $i<2$.
Since $\varrho$ is one-to-one we have 
  $K_0\cap K_1=\empt$, hence $\ell=\<K_0,K_1\>\in \Omega$, moreover
\begin{equation}\tag{$\star$}\label{eq:dis}
\Bigl(\bigcup_{m\in K_0}\dmem {\gamma}m\Bigr)
\cap
\Bigl(\bigcup_{m\in K_1}\dmem {\gamma}m\Bigr)=\empt
\end{equation}
because 
the elements of the family $\{B^{\gamma}_n(m):m<n\}$ are pairwise
disjoint. 

Now put $t=s\cup \{\<{\gamma},\ell\>\}$.
Then for each ${\varepsilon}\in {}^{\dom t}2$ we clearly have
\begin{equation}\tag{$\star\star$}\label{eq:diss}
 A_{{\varepsilon}({\gamma})}\cap \dcal[t,{\varepsilon}]\ne \empt,
\end{equation}
hence (\ref{eq:dis}) and (\ref{eq:diss}) together yield
that the extension $t$ of $s$ 
is $\dcal$-dyadic over $\blo{\gamma}=\blo {\ma(t)}$.
\end{proof}

Thus we may apply theorem \ref{tm:mc} to the matrix
$\dcal$ to obtain an infinite 
$\dcal$-dyadic partial function $h:{\kappa}\pto {\Omega}$.
Set $a=\dom h$ and  for each ${\alpha}\in a$ and $i<2$ 
put 
$L_{\alpha}^i=\{{\omega}{\alpha}+{\nu}(k):k\in \pii{h({\alpha})}\}$. 
For  $x\in L_{\alpha}^i$ put 
\begin{displaymath}
V(x)=\cup\{B^{\alpha}_{{\eta}(k)}({\nu}(k)):
x={\omega}{\alpha}+{\nu}(k)\text{ and }
k\in \pii{h({\alpha})}\}.
\end{displaymath}
Then $V(x)\in \bcal(x)$ because $\bcal(x)$ is closed under finite unions.
Since for $i<2$ 
\begin{equation}\notag
\bigl(\cup\{V(x):{x\in L^i_{\alpha}}\}\bigr)\cap {\kappa}
=\cup\{\dmem {\alpha}k:{k\in \pii {h({\alpha})}}\}
\end{equation}
and
\begin{displaymath}
\cup\{\dmem {\alpha}k:{k\in \pin {h({\alpha})}}\}\cap
\cup\{\dmem {\alpha}k:{k\in \pie {h({\alpha})}}\}=\empt,
\end{displaymath}
we have
\begin{displaymath}
\bigl(\cup\{V(x):{x\in L^i_{\alpha}}\}\bigr)\cap \bigl(\cup\{V(x):{x\in L^i_{\alpha}}\}\bigr)=\empt
\end{displaymath}
because the latter intersection is an open set which does 
not intersect the dense set $\blo0\subs{\kappa}$.
Hence 
the infinite family
\begin{equation}\notag
\Bigl\{\bigl\langle
\bigcup_{x\in L^0_{\alpha}}V(x),\bigcup_{x\in L^1_{\alpha}}V(x)
\bigr\rangle:{\alpha}\in a\Bigr\}
\end{equation}
is indeed dyadic.
\end{proof}

\begin{proof}[Proof of theorem \ref{tm:mc}]
\prlabel{tm:mc}

The proof will be based on the following 
two lemmas,   \ref{lm:nice} and \ref{tm:coh}. For these we need some more
notation and a new and rather 
technical notion of extendibility for set matrices.

Given a set $A$ we set
\begin{displaymath}
\ffcal {A}=\{f\in\fn(A,A):
\text{$f$ is injective and $\dom(f)\cap \ran(f)=\empt$}\}.
\end{displaymath}
Each function $f\in \ffcal{A}$ can be extended in natural way to
a bijection $f^*:A\to A$ as follows: 
\begin{displaymath}
f^*(a)=
\left\{
\begin{array}{ll}
f(a)&\text{if $a\in \dom f$,}\\
f^{-1}(a)&\text{if $a\in \ran f$,}\\
a&\text{otherwise}.
\end{array}
\right.
\end{displaymath}

\begin{definition}\label{df:nice} If $S$ and $T$ are sets of
ordinals then  the matrix
$\acal\in\matr ST$
 is called {\em \nice\ } iff 
for each $f\in \ffcal {S}$
there are a family $\nfam f\subs \sscal S$
and a function $\kfun {f}:\nfam {f}\to \br S;\le {\omega};$
such that
\begin{enumerate}[(1)] 
\item \label{di}  
the pair $(f,s)$ is $\acal$-dyadic whenever
$f\in \ffcal{S}$ and $s\in \nfam{f}$,
\item \label{empt} $\empt\in \nfam{f}$ for each $f\in \ffcal {S}$,
\item \label{open}
for  $f,g\in \ffcal S$ and $s\in \nfam f$
if $f^*\restr \kfun {f}(s)=g^*\restr \kfun {f}(s)$
then $s\in \nfam {g}$.
\item \label{ext} for any $f\in \ffcal {S}$, $s\in \nfam f$
and ${\alpha}\in S\cap \ma (s)$ there is
$\ell\in {\Omega}$  such that 
 $s\cup\{\<{\alpha},\ell\>\}\in \nfam f$.
\end{enumerate}
\end{definition}

Clearly, this last condition (\ref{ext}) is what explains our terminology.

\begin{lemma}\label{lm:nice}
If  ${\kappa}>\oo$ is regular and $\acal\in \omatr {\kappa}$ is a 
\nice\ matrix then 
there is an infinite partial function 
$h:{\kappa}\pto {\omega}$ 
that  is ${\acal}$-dyadic .
\end{lemma}

\begin{proof}\prlabel{lm:nice}
By induction on $n\in {\omega}$ we will define functions
$h_0\subs h_1\subs\dots h_n\subs\dots$ from $\sscal{\kappa}$
such that $|h_n|=n$ and
for each ${\nu}\in {\kappa}$
 
\begin{enumerate}[(i)]
\item[$(*)^{n}_{\nu}$] \label{in}
there is $g\in \ffcal {\kappa}$ such that $\ma (g)>{\nu}$,
$\ran g=\dom h_n$ and $h_n\circ g\in \nfam {g}$.
\end{enumerate}

First observe that  $h_0=\empt$
satisfies our requirements because,  according to (\ref{empt}),
condition  $(*)^{0}_{\nu}$  holds trivially 
for each ${\nu}<{\kappa}$.
 
Next assume that the construction has been done and the induction
hypothesis has been established   for $n$.
For each ${\nu}<{\kappa} $ 
choose a function $g_{\nu}\in \ffcal {\kappa}$
witnessing $(*)^n_{{\nu}+\oo}$ and then
write $K_{\nu}=\kfun {g_{\nu}}(h_n\circ g_{\nu})$ and
pick ${\zeta}_{\nu}\in ({\nu},{\nu}+\oo)\setm K_{\nu}$.
Clearly the set 
\begin{displaymath}
L=\{{\xi}\in {\kappa}:|\{{\nu}<{\kappa}:{\xi}\notin K_{\nu}\}|<{\kappa}\}
\end{displaymath}
is countable and so we can pick
${\xi}_{n}\in {\kappa}\setm (L\cup\dom h_n)$;
then the set
\begin{displaymath}
J=\{{\nu}<{\kappa}:
{\xi}_{n}\notin  K_{\nu}\}
\end{displaymath}
is of size ${\kappa}$.

Now set $g'_{\nu}=g_{\nu}\cup \{\<{\zeta}_{\nu},{\xi}_{n}\>\}$
for every ${\nu}\in J$.
For every such ${\nu}$ then ${\zeta}_{\nu},{\xi}_n\notin K_{\nu}$ implies
$g_{\nu}{}^*\restr K_{\nu}=g'_{\nu}{}^*\rest K_{\nu}$,
hence $h_n\circ g_{\nu}\in \nfam {g'_{\nu}}$ by (\ref{open}). 
Since ${\zeta}_{\nu}<{\nu}+\oo<
\ma (g_{\nu})=
\ma( h_n\circ g_{\nu})$, we can now apply 
(\ref{ext}) to get 
$\ell^{\nu}\in {\Omega}$ such that 
$(h_n\circ g_{\nu})\cup \{\<{\zeta}_{\nu},\ell^{\nu}\>\}\in \nfam{g'_{\nu}}$.

We can then fix $\ell_n\in {\Omega}$ such that 
$J_n=\{{\nu}\in J:\ell^{\nu}=\ell_n\}$ is of size ${\kappa}$
and
let $h_{n+1}=h_n\cup \{\<{\xi}_n,\ell_n\>\}$.

If ${\nu}\in J_n $ then
$h_{n+1}\circ g'_{\nu}=
(h_n\circ g_{\nu})\cup \{\<{\zeta}_{\nu},\ell_n\>\}
\in  \nfam{g'_{\nu}}$ and $\ma(g'_{\nu})>{\nu}$, so
$g'_{\nu}$ witnesses $(*)^{n+1}_{{\nu}}$.
But $J_n$ is unbounded in ${\kappa}$,
hence the inductive step is completed.

By $(*)^n_0$, for each $n<{\omega}$ there is $g_n$ such that
$\dom h_n=\ran g_n$ and 
$h_n\circ g_n\in \nfam{g_n}$. Hence, by (\ref{di}), $(g_n,h_n\circ g_n)$ is
$\acal$-dyadic, and so  $h_n$ is $\acal$-dyadic according to
observation \ref{obse:circ}. 
Consequently  $h=\bigcup\{h_n:n<{\omega}\}$ is as required: 
it is $\acal$-dyadic and infinite. 
\end{proof}

Given any infinite set $I$ we denote by $\coh I$ 
the poset $\fn(I,2)$,
i.e. the standard notion of forcing that adds $|I|$ many Cohen reals.

\begin{lemma}\label{tm:coh}
Let ${\kappa}=(2^{\omega})^+$.  
Then for each ${\lambda}$ we have
\begin{multline}\notag
V^{\coh {\lambda}}\models\text{If $\dcal\in\matr {\kappa}{\kappa}$ is
both \good\ and 
m.d.-extendible then there is $I\in \br {\kappa};{\kappa};$  }
\\\text{such that $\dcal^*=\<\dmem {\alpha}n\cap {\omega}:
\<{\alpha},n\>\in I\times {\omega}\>$ is \nice.}
\end{multline}
\end{lemma}

\begin{proof}\prlabel{tm:coh}
Assume that 
\begin{displaymath}
1_{\coh {\lambda}}\force\text{$\dot\dcal\in\matr {\kappa}{\kappa}$ is
m.d.-extendible.}
\end{displaymath}
Let ${\theta}$ be a large enough regular cardinal and consider the structure
$ {\hcal}_{\theta} =
\<H_{\theta},\in,\triangleleft,{\kappa},{\lambda},\dot\dcal\>$,
where $H_{\theta}=
\bigl\{x:|\operatorname{TC}(x)|<{\theta}\bigr\}$
and $\triangleleft$ is a fixed well-ordering of
$H_{\theta}$.

Working in $V$, for each ${\alpha}<{\kappa}$ choose a countable 
elementary submodel $N_{\alpha}$ of ${\hcal}_{\theta}$ 
with ${\alpha}\in N_{\alpha}$.
Then
there is  $I\in \br {\kappa};{\kappa};$ such that the models 
$\{N_{\alpha}:{\alpha}\in I\}$ are not only pairwise isomorphic but, denoting
by $\siso {\alpha}{\beta}$ the unique isomorphism between
$N_{\alpha}$ and $N_{\beta}$, we have
\begin{enumerate}[(i)]
\item\label{delta} 
the family
$\{N_{\alpha}\cap \theta:{\alpha}\in I\}$ forms a $\Delta$-system
with kernel $\Lambda$,
\item $\siso {\alpha}{\beta}({\xi})={\xi}$ for each ${\xi}\in \Lambda$,
\item $\siso{\alpha}{\beta}({\alpha})={\beta}$.
\end{enumerate}

For each ${\alpha}<{\kappa}$ and $n<{\omega}$ let
$\ddmem {\alpha}{n}$ be the ${\triangleleft}$-minimal
$\coh {\lambda}$-name of the $\<{\alpha},n\>^{\rm th}$ entry
of $\dot\dcal$.
Since ${\triangleleft} $ is in $\hcal_{\theta}$ and 
$\siso {\alpha}{\beta}({\alpha})={\beta}$ we have 
\begin{claim}
$\siso {\alpha}{\beta}(\ddmem {\alpha}n)=\ddmem {\beta}n$
for each ${\alpha},{\beta}\in I$ and $n\in {\omega}$.
\end{claim}

Let $G$ be any $\coh{\lambda}$-generic filter over $V$.
We shall show that 
\begin{displaymath}
V[G]\models\text{``$\dcal^*=\<\dmem {\alpha}n\cap {\omega}:
\<{\alpha},n\>\in I\times {\omega}\>$ is \nice.''}
\end{displaymath}

For each $f\in \ffcal {I}$ define the bijection
${\rho}_f:{\lambda}\to {\lambda}$ as follows:
\begin{displaymath}
{\rho}_f({\xi})=
\left\{
\begin{array}{ll}
\siso{\alpha}{f^*({\alpha})}({\xi})&\text{if ${\xi}\in N_{\alpha}\cap {\lambda}$
for some ${\alpha}\in I$,}\\
{\xi}&\text{otherwise.}
\end{array}
\right.
\end{displaymath}

In a natural way ${\rho}_f$ extends to  an automorphism of
$\coh {\lambda}$, 
 which will be denoted by ${\rho}_f$ as well.
Clearly, we have
\begin{claim}
If $f\in \ffcal I$, $f({\alpha})={\beta}$, $p\in \coh {\lambda}\cap
N_{\alpha}$ then $\siso {\alpha}{\beta}(p)={\rho}_f(p)$.
\end{claim}

For $f\in\ffcal {I}$ let $G^f=\{{\rho}^{-1}_f(p):p\in G\}$
and then set
\begin{multline}\notag
\nfam f=\{s \in \fn(I,\Omega) :\text{$s$ is $ \dot\dcal[G^f]$-min-dyadic}\}=
\\\{s \in \fn(I,\Omega) :\exists q\in G^f\ q\force 
\text{``$s$ is $ \dot\dcal$-min-dyadic''}\}.
\end{multline}

To define $K^f$, for each $s\in \nfam f$ pick a condition
$p_s\in G$ such that 
\begin{displaymath}
{\rho}^{-1}_f(p_s)\force \text{$s$ is $\dot\dcal$-min-dyadic}
\end{displaymath}
and let 
\begin{displaymath}
K^f(s)=\{{\alpha}\in I: (N_{\alpha}\setm \Lambda)\cap \dom p_s\ne \empt\}.
\end{displaymath}

Note that $K^f(s)$ as defined above is finite, although  
\ref{df:nice}.(\ref{open}) only requires 
$K^f(s)$ to be countable.

To check   property \ref{df:nice}.(\ref{open}) 
assume that $f,g\in \ffcal I$ and $s\in \nfam f$
with $g^*\restr K^f(s)=f^*\restr K^f(s)$.
Then ${\rho}^{-1}_g(p_s)={\rho}^{-1}_f(p_s)$
and so 
\begin{displaymath}
{\rho}^{-1}_g(p_s)\force\text{$s$ is $\dot\dcal$-min-dyadic},
\end{displaymath}
hence $s$ is also $\dot\dcal[G^g]$-min-dyadic , i.e.   $s\in \nfam g$.

Before checking \ref{df:nice}.(\ref{di}
)
we need one more observation.
\begin{claim}\label{lm:csere}
$\ddmem{f({\alpha})}{n}[G]\cap {\omega}=
\ddmem{\alpha}{n}[G^f]\cap {\omega}$ whenever $f\in \ffcal I$, 
${\alpha}\in \dom f$, and $n<{\omega}$.
\end{claim}

\begin{proof}[Proof of claim \ref{lm:csere}]
Let $k\in {\omega}$. Then 
$k\in \ddmem {f({\alpha})}{n}[G]$ iff
$\exists p\in G$ $p\force$ ``$k\in \ddmem {f({\alpha})}{n}$''
iff
$\exists p\in G\cap N_{f({\alpha})}$ $p\force$ ``$k\in \ddmem {f({\alpha})}{n}$''
iff
$\exists q\in G^f\cap N_{\alpha}$ $p=\siso{\alpha}{f({\alpha})}(q)\force$ ``$k\in \ddmem {f({\alpha})}{n}$''
iff $\exists q\in G^f\cap N_{\alpha}$ $q\force$ ``$k\in
\ddmem {\alpha}{n}$'' 
iff $\exists q\in G^f$ $q\force$ ``$k\in
\ddmem {\alpha}{n}$'' 
iff $k\in \ddmem{\alpha}{n}[G^f]$.
\end{proof}

Now let $f\in \ffcal I$ and $s\in \nfam f$.
By the definition of $\nfam f$, $s$ is $\dot\dcal[G^f]$-min-dyadic  
and so by observation \ref{lm:mdd}
$s$ is $\dot\dcal[G^f]$-dyadic over   ${\omega}$.
But it follows from \ref{lm:csere}, that
$s$ is $\dot\dcal[G^f]$-dyadic over   ${\omega}$ if and only if the pair 
$(f,s)$ is $\dot\dcal[G]$-dyadic over   ${\omega}$.

\ref{df:nice}.(\ref{empt}) is clear because 
$\empt$ is trivially $\acal$-min-dyadic for any $\acal\in \matr {\kappa}{\omega}$.
Finally \ref{df:nice}.(\ref{ext}) follows from the definition of $N(f)$
because
$\dot\dcal[G^f]$ is m.d.-extendible.
\end{proof}

Now, to complete the proof of theorem \ref{tm:mc}, first apply
 lemma \ref{tm:coh} to get 
$I\in \br {\kappa};{\kappa};$
such that 
\begin{displaymath}
\dcal^*=\<\dmem {\alpha}n\cap {\omega}:
\<{\alpha},n\>\in I\times {\omega}\>
\end{displaymath}
 is \nice. Then applying 
lemma \ref{lm:nice} to $\dcal^*$ we obtain
an infinite $\dcal^*$-dyadic function 
$h:{\kappa}\pto {\Omega}$. Since the matrix $\dcal$ is \good\
the function $h$  is $\dcal$-dyadic, as well.
\end{proof}

\section{Cardinal sequences of regular and $0$-dimensional spaces}
\label{sc:regular}

 For any 
regular, scattered  space $X$ we have $|X|\le 2^{|I(X)|}$, hence 
$\htt(X)<(2^{|I(X)|})^+$ 
and $|I_{\alpha}(X)| \le 2^{|\lev 0;X;|}$ for each $\alpha$.
This implies that for such a space $X$ its cardinal sequence $s$
satisfies $length(s) < (2^{|I(X)|})^+$ and $s(\alpha) \le 2^{s(\beta)}$
whenever $\beta < \alpha$. We shall show below that these properties of
a sequence $s$ actually characterize the cardinal sequences of regular
scattered spaces. 

In \cite{BS}, for each ${\gamma}<(2^{\omega})^+$, 
a 0-dimensional, scattered space of height ${\gamma}$ and
width ${\omega}$ was constructed. The next lemma generalizes that construction.

For an infinite cardinal ${\kappa}$, let $S_{\kappa}$ be the following
family of sequences of cardinals: 
\begin{displaymath}
S_{\kappa}=\bigl\{\<{\kappa}_i:i<{\delta}\>:{\delta}<(2^{\kappa})^+,\ 
{\kappa}_0={\kappa}\mbox{ and } {\kappa}\le {\kappa}_i\le
2^{\kappa}\mbox{ for each $i<{\delta}$ }
\bigr\}.
\end{displaymath}

\begin{lemma}\label{lm:const}
For any infinite cardinal ${\kappa}$ and $s\in S_{\kappa}$ there is 
$0$-dimensional scattered space $X$ with  $\operatorname{CS}(X)=s$.
\end{lemma}

\begin{proof}\prlabel{lm:const}
Let $s=\<{\kappa}_{\alpha}:{\alpha}<{\delta}\>\in S_{\kappa}$.
Write $X=\bigcup\bigr\{\{{\alpha}\}\times {\kappa}_{\alpha}:
{\alpha}<{\delta}\bigr\}$. 
Since $|I|\le 2^{\kappa}$ we can fix an independent family
$\{F_x:x\in X\}\subs \br {\kappa};{\kappa};$.

The underlying set of our space is $X$ and and the topology ${\tau}$ on $X$ is
given by declaring for each 
$x=\<{\alpha},{\xi}\>\in X$
the set 
\begin{displaymath}
U_x=\{x\}\cup ({\alpha}\times F_x) 
\end{displaymath}
to be clopen, i.e. $\{U_x,X\setm U_x:x\in X\}$ is a subbase for ${\tau}$. 

The space $X$ is clearly $0$-dimensional and $T_2$.

\begin{claim}
If $x=\<{\beta},{\xi}\>\in U\in {\tau}$ and ${\alpha}<{\beta}$ then
$U\cap (\{{\alpha}\}\times {\kappa}_{\alpha})$ is infinite.
\end{claim}

\begin{proof}[Proof of the claim]\prnolabel
We can find disjoint sets $A,B\in \br X\setm \{x\};<{\omega};$ such that  
\begin{displaymath}
x\in U_x\cap\bigcap_{y\in A}U_y\setm \bigcup_{z\in B}U_z\subs U.
\end{displaymath}
Observe that if $\<{\gamma},{\xi}\>\in A$ then
${\beta}<{\gamma}$. Thus
\begin{equation}\notag
U\cap (\{{\alpha}\}\times {\kappa}_{\alpha})\supset
\{{\alpha}\}\times\Bigl(\bigcap_{y\in A\cup\{x\}}F_y\setm \bigcup_{z\in B} F_z\Bigr),
\end{equation}
and the  set on the right side is infinite because $\{F_x:x\in X\}$ was chosen
to be independent.
\end{proof}

To complete our proof,  by induction on ${\alpha}<{\kappa}$, we verify 
that  $\lev{\alpha};X;= 
\{{\alpha}\}\times {\kappa}_{\alpha}$, hence $\operatorname{CS}(X)=s$.
Assume that this is true for ${\nu}<{\alpha}$. 
If $x\in \{{\alpha}\}\times {\kappa}_{\alpha}$ then
\begin{displaymath}
U_x\cap \bigl(X\setm \bigcup_{{\nu}<{\alpha}}\lev {\nu};X;\bigr)=\{x\},
\end{displaymath}
hence $\{{\alpha}\}\times {\kappa}_{\alpha}\subs \lev {\alpha};X;$.
On the other hand, if $x=\<{\beta},{\xi}\>\in X$ with  ${\beta}>{\alpha}$ and 
$U\in {\tau}$ is a neighbourhood of $x$, then, by the claim above, 
$U\cap (\{{\alpha}\}\times {\kappa}_{\alpha})$ is infinite, hence $x$ is
not
isolated in $X\setm \bigcup_{{\nu}<{\alpha}}\lev {\nu};X;$, i.e.,
$x\notin \lev {\alpha};X;$.
Thus $\lev {\alpha};X;=\{{\alpha}\}\times {\kappa}_{\alpha}$.
\end{proof}

\begin{theorem}\label{tm:reg} For any sequence $s$ of cardinals 
the following statements are equivalent
\begin{enumerate}[(1)]
\item \label{r}$s=CS(X)$ for some 
 regular scattered space $X$, 
\item \label{0}$s=CS(X)$ for some  
0-dimensional scattered space $X$, 
\item \label{c}
 for some natural number $m$ there are infinite cardinals
${\kappa}_0>{\kappa}_1>\dots>{\kappa}_{m-1}$ and for all $i<m$ sequences 
$s_i\in S_{{\kappa}_i}$ such that 
$s=s_0{}^\frown s_1{}^\frown \dots{}^\frown s_{m-1}$ or
$s=s_0{}^\frown s_1{}^\frown \dots{}^\frown s_{m-1}{}^\frown \<n\>$ 
for some natural number $n>0$.
\end{enumerate}
\end{theorem}

\begin{proof}\prlabel{tm:reg}\mbox{ } \\
(\ref{r})$\Longrightarrow$ (\ref{c})\\
By induction on $j$ we choose ordinals 
${\nu}_j<\htt(X)$ and cardinals ${\kappa}_j$ such that ${\nu}_0=0$
and ${\kappa}_0=|\lev 0;X;|$, moreover, for $j>0$ with
${\kappa}_{j-1}$ infinite 
\begin{displaymath}
{\nu}_j=\min\bigr\{{\nu}\le\htt(X):
|\lev {\nu};X;|<{\kappa}_{j-1}\},
\end{displaymath}
and   
${\kappa}_j=|\lev {\nu}_j;X;|$.
We stop when ${\kappa}_m$ is finite.
For each $j<m$ let ${\delta}_j={\nu}_{j+1}\stackrel{.}{-}{\nu_j}$.
Then the sequence $s_j=\<|\lev {{\nu}_j+{\delta}};{X};|:{\delta}<{\delta}_j\>$
is in $S_{{\kappa}_j}$.
Thus $\operatorname{CS}(X)=s_0{}^\frown s_1{}^\frown \dots{}^\frown s_{m-1}$ provided 
${\kappa}_m=0$ (i.e. $\lev {\nu}_m;X;=\empt$)  and 
$\operatorname{CS}(X)=s_0{}^\frown s_1{}^\frown \dots{}^\frown s_{m-1}{}^\frown \<{\kappa}_m\>$
when $0<{\kappa}_m<{\omega}$.

\medskip

\noindent
(\ref{c})$\Longrightarrow$ (\ref{0})\\
First we prove this implication  for sequences $s$ of the form 
$s_0{}^\frown s_1{}^\frown \dots{}^\frown s_{m-1}$ 
by induction on $m$.
If  $s\in S_{{\kappa}_o}$ then 
the statement is just lemma \ref{lm:const}

Assume now that 
$s=s_0{}^\frown s_1{}^\frown \dots{}^\frown s_{m-1}$, where 
 ${\kappa}_0>{\kappa}_1>\dots>{\kappa}_{m-1}$ and 
$s_i\in S_{{\kappa}_i}$ for $i<m$.

According to lemma \ref{lm:const} there is a $0$-dimensional space $Y$  
  with cardinal sequence $s_{m-1}$.
Using the  inductive assumption we can also fix pairwise disjoint
$0$-dimensional topological spaces $X_{y,n}$ for  
$\<y,n\>\in I_0(Y)\times {\omega}$, each having the cardinal 
sequence $s'=s_0{}^\frown s_1{}^\frown \dots{}^\frown s_{m-2}$.
We then define the  space $Z=\<Z,{\tau}\>$ as follows. Let 
\begin{displaymath}
Z=Y\cup\bigcup\{X_{y,n}: y\in I_0(Y), n<{\omega}\}.
\end{displaymath}
A set $U\subs Z$ is in ${\tau}$ iff
\begin{enumerate}[(i)]
\item $U\cap Y$ is open in $Y$, 
\item $U\cap X_{y,n}$ is open in $X_{y,n}$ for each 
$\<y,n\>\in I_0(Y)\times {\omega}$,
\item if $y\in I_0(Y)\cap U$ then there is $m<{\omega}$ such that 
$\bigcup\{X_{y,n}:m<n<{\omega}\}\subs U$.
\end{enumerate}
If $U$ is a clopen subset of $Y$ and $n<{\omega}$ then it is easy to
check that 
\begin{equation}\notag
Z(U,n)=U\cup\bigcup\{X_{y,m}:y\in \lev 0;Y;\cap U, n<m<{\omega}\}
\end{equation}
is clopen in $Z$. Hence 
\begin{multline}\notag
\bcal=\{Z(U,n):U\subs Y\ \text{ is clopen}, n<{\omega}\}\cup\\
\{T:\text{$T$ is a clopen subset of some $X_{y,n}$}\}
\end{multline}
is a clopen base of $Z$ and so
$Z$ is $0$-dimensional.

Let ${\delta}'=\operatorname{length}(s')$ and 
${\delta}=\operatorname{length}(s)$.
\begin{claim}\label{cl:lev1}
$\lev {\alpha};Z;=\bigcup\{\lev {\alpha};X_{y,n};:\<y,n\>\in \lev
0;Y;\times {\omega}\}$ for ${\alpha}<{\delta}'$.
\end{claim}
\begin{proof}[Proof of the claim \ref{cl:lev1}]\prlabel{cl:lev1}
Since $X_{y,n}$ is an open subspace of $Z$ it follows that
$\lev {\alpha};X_{y,n};\subs \lev {\alpha};Z;$. On the other hand,
\begin{equation}\notag
Y\subs{\overline{\bigcup\{\lev {\alpha};X_{y,n};:\<y,n\>\in \lev
0;Y;\times {\omega}\}}}^Z,
\end{equation}
 hence $Y\cap \lev{\alpha};Z;=\empt$. 
\end{proof}
Since, by  claim \ref{cl:lev1},
\begin{equation}\notag 
Z\setm \bigcup_{{\alpha}<{\delta}'}\lev {\alpha};Z;=Y,
\end{equation}
it follows that for ${\delta}'\le {\alpha}<{\delta}$ we have
\begin{equation}\tag{$*$}\label{eq:lev2}
\lev {\alpha};Z;=\lev {\alpha}\stackrel{.}{-}{\delta}';Y;.
\end{equation}

Thus $Z=\bigcup_{{\alpha}<{\delta}}\lev {\alpha};Z;$, hence 
$Z$ is a scattered space of height ${\delta}$.

If ${\alpha}<{\delta}'$ then, by  claim \ref{cl:lev1}, 
\begin{displaymath}
|\lev {\alpha};Z;|=|\lev 0;y;|\cdot {\omega}\cdot s'({\alpha})=
{\kappa}_{m-1}\cdot {\omega}\cdot 
s'({\alpha})=s'({\alpha})=s({\alpha}).
\end{displaymath}
If ${\delta}'\le {\alpha}<{\delta}$ then, by  (\ref{eq:lev2}),  
$|\lev {\alpha};Z;|=|\lev {\alpha}\stackrel{.}{-}{\delta}';Y;|=
s_{m-1}({\alpha}\stackrel{.}{-}{\delta}')=s({\alpha})$, consequently
$\operatorname{CS}(Z)=s$.

Thus  we proved the statement for sequences of the form 
$s_0{}
^\frown \dots{}^\frown s_{m-1}$.

If $s=s_0{}^\frown \dots{}^\frown s_{m-1}\,{}^\frown \<n\>$ then
writing $s'=s_0{}^\frown \dots{}^\frown s_{m-1}$ we can first
find pairwise disjoint $0$-dimensional scattered spaces $X_{i,m}$, 
$\<i,m\>\in n\times {\omega}$ each having cardinal 
sequence $s'$. 
Let 
\begin{displaymath}
Z=\{x_i:i<n\}\cup\bigcup\{X_{\<i,m\>}: i<n, m<{\omega}\}.
\end{displaymath}
Declare a set $U\subs Z$  open iff
\begin{enumerate}[(i)]
\item $U\cap X_{i,m}$ is open in $X_{i,m}$ for each 
$\<i,m\>\in n\times {\omega}$,
\item if $x_i\in U$ then there is $n_i<{\omega}$ such that 
$\bigcup\{X_{i,m}:n_i<m<{\omega}\}\subs U$.
\end{enumerate}
Then $Z$ is 0-dimensional, and
\begin{equation}\notag
\lev{\alpha};Z;=\left\{
\begin{array}{ll}
\bigcup\{\lev {\alpha};X_{i,m};:i<n, m<{\omega}\}&
\text{if }{\alpha}<\operatorname{length}(s'),\\
\{x_i:i<n\}&\text{if }{\alpha}=\operatorname{length}(s').
\end{array}
\right.
\end{equation}
Hence again $Z$ is  a scattered space with 
$\operatorname{CS}(Z)=s$.

\medskip

\noindent
(\ref{0})$\Longrightarrow$ (\ref{r})
Straightforward.

\end{proof}

We leave it to the reader to verify that the sequences described in item
(3) of theorem \ref{tm:reg} are exactly those mentioned in the beginning
of the section with the additional obvious necessary condition that all
but the last term of the sequence are infinite cardinals.

\end{document}